\documentclass[11pt,a4paper]{article}
\usepackage{amssymb,comment,amsmath,amsthm}
\usepackage{graphicx}
\usepackage{color}
\usepackage{titling,titlesec}
\usepackage{url}
\usepackage{amsthm,lipsum}
\usepackage[nocompress]{cite}

\newfont{\blb}{msbm10 scaled\magstep1}
\newfont{\comp}{cmr12 scaled\magstep1}
\newfont{\compb}{cmr10 scaled\magstep2}
\newfont{\sbb}{cmssbx10 scaled\magstep3}
\newfont{\sbbb}{cmssbx10 scaled\magstep5}
\newfont{\sbs}{cmssbx10 scaled\magstep1}
\newtheorem{theorem}{Theorem}

\newtheorem{conj}{Conjecture}

\newcommand{\bs}{{\scriptstyle \blacksquare}}

\title{Independent sets in hypergraphs}

\author{
Jacques Verstraete\thanks{Department of Mathematics, University of California, San Diego, CA, 92093-0112 USA.
			Research supported by NSF awards DMS-1952786 and DMS-2347832. 
E-mail: jacques@ucsd.edu.} \quad \quad Chase Wilson\footnotemark[1] \thanks{E-mail: c7wilson@ucsd.edu.}}

\begin{document}

\date{ }

\maketitle

\begin{abstract}
A theorem of Shearer states that every $n$-vertex triangle-free graph of maximum degree $d \geq 2$ contains an independent set 
of size at least $(d\log d - d + 1)/(d - 1)^2 \cdot n$. 
Ajtai, Koml\'{o}s, Pintz, Spencer and Szemer\'{e}di proved that every $(r + 1)$-uniform  $n$-vertex ``uncrowded'' hypergraph of maximum degree $d \geq 1$
has an independent set of size at least $c_r(\log d)^{1/r}/d^{1/r} \cdot n$ for some $c_r > 0$ depending only on $r$. Shearer asked whether his method for triangle-free graphs could be extended to uniform hypergraphs. 
In this paper, we answer this in the affirmative, thereby giving a short proof of the theorem of Ajtai, Koml\'{o}s, Pintz, Spencer and Szemer\'{e}di for a wider class of 
``locally sparse'' hypergraphs.
\end{abstract}

\section{Introduction}

The {\em independence number} $\alpha(H)$ of a hypergraph $H$
is the maximum size of an {\em independent  set} of vertices of $H$ -- a set of vertices containing no edges of $H$. For $k \geq 2$, a {\em $k$-cycle} is a hypergraph consisting of edges $e_0,e_1,\dots,e_{k-1}$ for which there exist distinct vertices $v_0,v_1,\dots,v_{k-1}$ such that $v_i \in e_i \cap e_{i + 1}$ for $0 \leq i < k$ with subscripts mod $k$. A hypergraph $H$ is {\em linear} if it contains no 2-cycle -- equivalently, for any edges $e,f$ in $H$, $|e \cap f| \leq 1$. 
A hypergraph is {\em locally sparse} if it contains no 2-cycle or 3-cycle. Following~\cite{AKPSS}, a hypergraph is {\em uncrowded} if it has no $k$-cycles for $k \leq 4$. 
Independent sets in linear, locally sparse and uncrowded hypergraphs have been studied at length and have numerous connections and applications, for instance to independent sets in Steiner systems~\cite{deBR,DLR,PR,RS,KMV}, Heilbronn's triangle problem~\cite{CPZ,KPS,R}, line arrangements~\cite{Fu,FuP}, Ramsey theory~\cite{A,AKS,BK,BPM,K,KMV2},
coloring hypergraphs~\cite{AKSu,BFM,CM,FM,FM2} and the method of containers~\cite{BMS,Samotij,ST}. 

\medskip

A theorem of Shearer~\cite{Sh} states that every $n$-vertex triangle-free (locally sparse) graph of average degree $d \geq 2$ contains an independent set 
of size at least $(d\log d - d + 1)/(d - 1)^2 \cdot n$, improving an earlier result of Ajtai, Koml\'{o}s and Szemer\'{e}di~\cite{AKS}. 
For uniform hypergraphs, a fundamental theorem due to Ajtai, Koml\'{o}s, Pintz, Spencer and Szemer\'{e}di~\cite{AKPSS} states
that for any $n$-vertex $(r + 1)$-uniform uncrowded hypergraph $H$ of maximum degree $d$, as $d \rightarrow \infty$:
\begin{equation}\label{basic2}
\alpha(H) \geq \frac{0.98}{e} \cdot 10^{-\frac{5}{r}} \cdot \Bigl(\frac{\log d}{d}\Bigr)^{\frac{1}{r}} \cdot n.
\end{equation}
The proof of (\ref{basic2}) is via a randomized greedy algorithm akin to the
R\"{o}dl nibble~\cite{R}, although the analysis is substantially more complicated. This result improves the bound~\cite{CT} obtained by random sampling
in any $n$-vertex $(r + 1)$-uniform hypergraph of maximum degree $d \geq 1$, namely:
\begin{equation}\label{basic1}
\alpha(H) \geq \Bigl(1 - \frac{1}{r}\Bigr)\frac{n}{d^{1/r}}.
\end{equation}
Shearer asked whether the method in~\cite{Sh} for triangle-free (locally sparse) graphs could apply to the uniform hypergraph setting. We answer this in the affirmative, 
and thereby give a short proof of the theorem Ajtai, Koml\'{o}s, Pintz, Spencer and Szemer\'{e}di~\cite{AKPSS} for the slightly wider class of locally sparse hypergraphs:

\begin{theorem} \label{main}
For $r \geq 1$, there exists $c_r > 0$ such that for any $(r + 1)$-uniform $n$-vertex locally sparse hypergraph $H$ with maximum degree $d \geq 1$. 
\begin{equation}\label{mainbound}
	\alpha(H) \geq  c_r \cdot \Bigl(\frac{\log d}{d}\Bigr)^{\frac{1}{r}} \cdot n.
 \end{equation}
where $c_r = (1 - o_d(1)) f(r)$ for some function $f$.
\end{theorem}

The proof of Theorem \ref{main} is inspired by the approach of Kostochka, Mubayi and the first author~\cite{KMV} and  Shearer~\cite{Sh} for triangle-free graphs, and gives an improvement of the constant 
factor in~\cite{AKPSS} for small $r$, for instance, it is possible from the proof of Theorem \ref{main} to obtain $c_2 = 1/8$ for large enough $d$, which improves on the constant value $0.98/10^{5/2}e \approx 0.00114007$ in (\ref{basic2}). The value of $c_r$ from the proof of Theorem \ref{main} is 
\begin{equation} 
c_r = (1 - o_d(1)) \frac{(r - 1)^{1 - \frac{1}{r}} (2^r - 1)}{r^{2 + \frac{2}{r}} 2^r} 
\end{equation}
and, in particular, $c_r = (1 - o_r(1))(1 - o_d(1))/r$ as $r,d \rightarrow \infty$. Via random sampling, as in Duke, Lefmann and R\"{o}dl~\cite{DLR}, one may obtain 
a lower bound which is a factor $(1 - o_d(1))(1 - 2/(3r - 1))^{1/r}$ times the lower bound (\ref{mainbound}) in Theorem \ref{main} when $H$ is a linear $(r + 1)$-uniform hypergraph of maximum degree $d$. In the case of $n$-vertex $(r + 1)$-uniform $d$-regular hypergraphs $H$ with no $k$-cycles for $k \leq 2g + 1$, it was shown by Nie and the first author~\cite{NV} that the randomized greedy algorithm for independent sets terminates with an independent set $I \in \mathcal{I}(H)$ such that as $d \rightarrow \infty$:
\begin{equation} 
(1 - o_d(1)) \cdot \Bigl(\frac{\log d}{rd}\Bigr)^{\frac{1}{r}} \cdot n - \frac{d^{g}}{rg!} \cdot n \leq \mathrm{E}(|I|) \leq (1 + o_d(1)) \cdot \Bigl(\frac{\log d}{rd}\Bigr)^{\frac{1}{r}} \cdot n + \frac{d^{g}}{rg!} \cdot n.
\end{equation}
In particular, if $g \geq d$, this gives an independent set of size $(1 + o_d(1)) (\log d)^{1/r}/(rd)^{1/r} \cdot n$.

\bigskip\medskip

{\bf Notation and terminology.} Let $H$ be a $(r + 1)$-uniform hypergraph with vertex set $V = V(H)$ and edge set $E(H)$.  We write
$\mathcal{I}(H)$ for the set of independent sets in $H$, and let $\alpha(H) = \max\{|I| : I \in \mathcal{I}(H)\}$ denote the {\em independence number} of $H$.
For $v \in V$, and $1 \leq k \leq r$, define
\[ N_k(v) = \{f \in { V \setminus \{v\}  \choose k } : \exists e \in E(H), f \cup \{v\} \subseteq e\}.\]
In particular, we let $N(v) = N_1(v)$ denote the set of vertices in a common edge with $v$ -- sometimes called the {\em neighborhood} of $v$. Define for $I \subseteq V$ the {\em $k$-shadow of $H$ on $I$} as
\[
\partial_k I = \{e \in {I \choose k} : \exists f \in E(H), e \subseteq f\} = \displaystyle{\bigcup_{v \in V}} N_k(v) \cap {I \choose k}.
\]
We use the following asymptotic notation. If for each $r \in \mathbb N$, $f_r,g_r : \mathbb N \rightarrow \mathbb R^+$ are positive real-valued functions, then we write $f_r = O_r(g_r)$ and $g_r = \Omega_r(f_r)$ if there exists for each $r \in \mathbb N$
a constant $C_r$ such that $f_r(x) \leq C_r g_r(x)$ for all $x \in \mathbb N$, and
$f_r = o_r(g_r)$ if $\lim_{r \rightarrow \infty} C_r = 0$.

\section{A Non-Uniform Distribution on Independent Sets}\label{intuition}

This section serves as intuition and motivation for the proof of Theorem \ref{main}. If $H$ is a locally sparse $n$-vertex $(r + 1)$-uniform hypergraph of maximum degree $d$, 
and we sample an independent set $I \in \mathcal{I}(H)$ uniformly, then following Shearer~\cite{Sh} (see page 272 in Alon and Spencer~\cite{AS}), it is natural to consider for $v \in V(H)$ the random variable $X_v$ which depends on the $r$-shadow $\partial_r I$ of $H$ on $I$:
\begin{equation}\label{def:xv} 
X_v = \left\{\begin{array}{ll}
d^{1/r} & \mbox{ if }v \in I \\
|N_r(v) \cap \partial_r I| & \mbox{ if }v \not \in I.
\end{array}\right.
\end{equation}
Let $X = \sum_{v \in V(G)} X_v$. By definition, $X = d^{1/r} \cdot |I| + |\partial_r I|$. To derive a lower bound on $\mathrm{E}(X_v)$, let 
\[ I_v = I \backslash (N(v) \cup \{v\}) \quad \mbox{ and } \quad Y_v = |\{e \in N_r(v) : I_v \cup e \in \mathcal{I}(H)\}|.\]
Since $H$ is linear, $N_r(v)$ consists of disjoint $r$-element sets, and since $H$ is locally sparse, if $I_v \cup e$ and $I_v \cup f$ are independent sets in 
$H$ for $e,f \in N_r(v)$, then so is $I_v \cup e \cup f$. We conclude
\[ \mathrm{E}(X_v | I_v) = \frac{Y_v 2^{r(Y_v - 1)} + d^{1/r} (2^r - 1)^{Y_v}}{2^{rY_v} + (2^r - 1)^{Y_v}}.\]
A simple optimization shows this is always at least $b_r\log d$ for some positive constant $b_r = (1 - o_r(1))/r$ depending only on $r$. By the tower property of conditional expectation, $\mathrm{E}(X_v) = \mathrm{E}(\mathrm{E}(X_v|I_v))$. Therefore, if $X = \sum_{v \in V(H)} X_v$, 
\[ b_r \cdot n\log d \leq  \mathrm{E}(X) \leq d^{1/r} \cdot \mathrm{E}(|I|) + \mathrm{E}(|\partial_r I|).   \]
Now $|\partial_r I|$ may be as large as $d|I|$, and if we use this upper bound, we only obtain the very weak bound 
$\mathrm{E}(|I|) \geq b_r(\log d)/(d + d^{1/r}) \cdot n$. The value of $\mathrm{E}(|\partial_r I|)$ seems to be difficult 
to determine when $I$ is a uniformly chosen independent set in a linear or locally sparse hypergraph. 

\medskip

To remedy the lack of control of $\mathrm{E}(|\partial I|)$, we shall sample independent sets in a locally sparse $(r + 1)$-uniform hypergraph with a carefully chosen 
non-uniform probability measure $\mathrm{Pr}(I)$, which weights our sampling of $I \in \mathcal I(H)$ by $\exp(-\delta|\partial_r I|)$, where $\delta^r \approx (\log d)/d$. 
This has the effect that the probability that $\partial_r I$ is large is vanishingly small. But this introduces a new problem: we may have that $N_r(v) \cap \partial_r I$ is empty but $N_{r - 1} \cap \partial_{r - 1} I$ is large with high probability. In this case, we do not get anything from the $|N_r(v) \cap \partial_r I|$ term in $X_v$, but also adding $v$ to $I$ decreases its weight by enough so that $\Pr(v \in I)$ is small. To fix this, we add a term counting $|N_{r - 1}(v) \cap \partial_{r - 1} I|$ to $X_v$. This forces us to weight our sampling of $I \in \mathcal I(H)$ by $\exp(-2\delta^2|\partial_{r - 1} I|)$ in order to control the upper bound on $X_v$. This effect cascades and leads us to the above definitions for $\Pr(I)$ and $X_v$ which depend on all of $|\partial_k I|$ where $1 \leq k \leq r$ -- see definitions (\ref{measure}) and (\ref{xv}) in the next section. The analysis for the lower bound then becomes, either there is some $k$ for which $|N_k(v) \cap \partial_k I|$ is large in which case we achieve the desired lower bound from the corresponding term, or $|N_k(v) \cap \partial_k I|$ is small for all $k$ in which case $v$ can be added to $I$ with good probability. This gives our intuition for the proof of Theorem \ref{main}, which we now present.

\section{Proof of Theorem \ref{main}}

\subsection{Definition of parameters}

Let $H$ be a $(r + 1)$-uniform locally sparse hypergraph of maximum degree $d \geq 1$ with $\alpha(H) = m$. For $2 \leq k \leq r$, set the following parameters:
\begin{equation}\label{parameters}
\alpha = d^{-\frac{1}{r}} (\log d)^{\frac{1}{r}}, \quad \qquad \delta^r = \frac{\alpha^r}{r^2(r - 1)}, \quad \qquad \beta_k = \exp(-(r - k + 1)\delta^{r-k+1})
\end{equation}
and let $\beta_0 = \beta_1 = 1$ and $\beta_{r + 1} = 0$.
We sample each $I \in \mathcal{I}(H)$ with probability
\begin{equation}\label{measure}
	\mathrm{Pr}(I) = \frac{\prod_{k = 2}^{r} \beta_k^{|\partial_{k} I|}}{\sum_{J \in \mathcal{I}(H)} \prod_{k = 2}^{r} \beta_k^{|\partial_{k } J|}}.
\end{equation}
For $I \in \mathcal{I}(H)$, let $X = X(I) = \sum_{v \in V(H)} X_v(I)$ where $X_v = X_v(I)$ is defined by
\begin{equation} \label{xv}
	X_v = \begin{cases}
		d^{\frac{1}{r}} & v \in I\\
		\sum_{k = 1}^{r} \delta^{r-k} |N_k(v) \cap \partial_k I| & v \not \in I.
	\end{cases}
\end{equation}

To prove Theorem \ref{main}, we find an upper and lower bound for $\mathrm{E}(X)$. We use the careful definition of the probability measure
$\mathrm{Pr}(I)$ to show $\mathrm{E}(X) \leq (1 + o_d(1))\alpha(H)(\log d)/\delta(r - 1)$ -- this is Claim 2. On the other hand, we show for each $v \in V(H)$,
$\mathrm{E}(X_v) \geq (1 - o_d(1))(2^r - 1)\log d/r^2 2^r$ -- this is Claim 5. Combining these bounds gives Theorem \ref{main}.

\subsection{The upper bound on $\mathrm{E}(X)$}

In this section, we prove $\mathrm{E}(X)  \leq (1 + o_d(1))(m\log d)/\delta(r - 1)$ as $d \rightarrow \infty$.

\medskip
\bigskip

{\bf Claim 1.} Let $Z = \log \prod_{k = 2}^r \beta_k^{-|\partial_k I|}$. Then
\begin{equation}
X \leq \Bigl(d^{\frac{1}{r}} + \frac{\log d}{\delta r(r - 1)}\Bigr) |I| + \frac{Z}{\delta}.
\end{equation}

{\it Proof.} Since $|N_1(w)| \leq dr$ for $w \in V(H)$, $\delta^{r-1}|N_1(w)| \leq \delta^{r-1}rd$.
As $H$ is linear, for $k \geq 2$ each element of $\partial_k I$ is in $N_k(w)$ for at most $r - k + 1$ vertices $w \in V(H) \backslash I$.
By definition (\ref{xv}),
\begin{eqnarray*}
	X & \leq & d^{\frac{1}{r}}|I| +  \sum_{k = 1}^r \sum_{v \in V} \delta^{r-k} |N_k(v) \cap \partial_k I| \\
&\leq& d^{\frac{1}{r}}|I| + r\delta^{r - 1}d|I| + \frac{1}{\delta}\sum_{k = 2}^r (r - k + 1)\delta^{r-k+1} |\partial_k I|.
\end{eqnarray*}
By definition (\ref{parameters}) of $\beta_k$, the last term is exactly $Z/\delta$. Since $|I| \leq \alpha(H)$, using
the definition (\ref{parameters}) of $\delta$, we have $r\delta^{r - 1}d = (\log d)/\delta r(r - 1)$,
which proves Claim 1. $\bs$

\bigskip
\medskip

{\bf Claim 2.} As $d \rightarrow \infty$,
\begin{equation}
\mathrm{E}(X)  \leq (1 + o_d(1)) \cdot \frac{m\log d}{\delta(r - 1)}.
\end{equation}

\medskip

{\it Proof.}  By definition of $\beta_k$ and (\ref{measure}), if $Z(I) \geq m/r \log d + 4m$ then 
\[
\Pr(I) =   \frac{ \exp ({ - Z(I)}) }{\sum_{J \in \mathcal I(H) } \prod_{k = 2}^r \beta_k^{|\partial_k J|}} \leq \frac{ \exp ( -m/r \log d - 4m) } { \prod_{k = 2}^r \beta_k^{|\partial_k \emptyset|}  } = \exp(   -m/r \log d - 4m)
\]
 Since $|\mathcal{I}(H)| \leq 2^m {n \choose m}$,
\begin{eqnarray*}
\Pr\Bigl(Z(I) \geq \frac{m}{r} \log d + 4m\Bigr) &\leq&
\exp\Bigl(-\frac{m}{r} \log d - 4m\Bigr) \cdot 2^m {n \choose m}.
\end{eqnarray*}
Using ${n \choose m} \leq (en/m)^m$ and the basic lower bound $m \geq (1 - 1/r)n/d^{1/r}$ from (\ref{basic1}),
\begin{eqnarray*} 2^m {n \choose m} &\leq&
 \exp\Bigl(m\log 2 + m\log\frac{en}{m}\Bigr) \\
 &\leq& \exp\Bigl(\frac{m}{r}\log d + m\log 2 + m + m\log\frac{r}{r - 1}\Bigr) \; \; \leq \; \; \exp\Bigl(\frac{m}{r}\log d + 3m\Bigr).
 \end{eqnarray*}
Therefore the above probability is at most $\exp(-m) = o(m^{-r})$ as $d \rightarrow \infty$.
Since $Z = O(m^r)$, we conclude $\mathrm{E}(Z) \leq (1 + o_d(1))m(\log d)/r$. Using Claim 1 and $d^{1/r} = o( \log d / (r \delta))$,
\[ \mathrm{E}(X) \leq d^{\frac{1}{r}}m + \frac{m\log d}{\delta r(r - 1)} + (1 + o_d(1))\frac{m\log d}{\delta r} = (1 + o_d(1)) \cdot \frac{m\log d}{\delta (r - 1)}.\]
This proves Claim 2. $\bs$

\subsection{The lower bound on $\mathrm{E}(X)$}

In this section we give a lower bound on $\mathrm{E}(X_v)$ for each $v \in V(H)$. Fixing a vertex $v \in V(H)$, let
\[ I_v  = I \backslash (N(v) \cup \{v\}).\]
For an edge $e \in E(H)$ with $v \in e$, let 
\[ \omega_e = \{x \in e \backslash \{v\} : I_v \cup \{x\} \in \mathcal{I}(H)\} \quad \mbox{ and } \quad Y_v = \{\omega_e : v \in e \in E(H)\}.\] 
For $\omega_e \in Y_v$ and  $S \subseteq \omega_e$, let
\begin{eqnarray*}
\beta_S =  \prod_{x \in S} \prod_{k = 1}^{r-1}\beta_{k + 1}^{|N_k(x) \cap I_v|} \qquad \qquad f_S = \prod_{j = 0}^{|S|} \beta_j^{{|S| \choose j}} \qquad \qquad f_S^+ = \prod_{j = 0}^{|S|} \beta_{j + 1}^{{|S| \choose j}}.
\end{eqnarray*}
The quantity $f_S$ corresponds to the weight in (\ref{measure}) coming from the edge $e$ after adding $S$ to $I_v$ and $f_S^+$ corresponds to the weight after adding $S \cup \{v\}$ to $I_v$. $\beta_S$ corresponds to the weight in (\ref{measure}) coming from edges of the form $ x \cup T$ where $x \in S$ and $T \subseteq I_v$ after adding $S$ to $I_v$. Since $H$ is locally sparse, this is the only additional weight in ({\ref{xv}) we will incur from adding $S$ or $S \cup \{v\}$ to $I_v$. In light of this, for $\omega \in Y_v$, define:
\begin{eqnarray*}
W _\omega = \sum_{S \subseteq \omega} \beta_S f_S \qquad \qquad W^+_\omega = \sum_{S \subseteq \omega} f_S^+ \beta_S
\end{eqnarray*}
$W_\omega$ is the total weight incurred when we do not add $v$ to $I$ and $W_\omega^+$ is the total weight incurred when we do add $v$ to $I$.

\medskip

{\bf Claim 3.} For $v \in V(H)$ and $I \in \mathcal{I}(H)$, let  $I_v = I \backslash \bigcup N_1(v) \cup \{v\}$
and let $A_v$ be the event $v \in I$ conditional on $I_v$ and let $A_v^c$ be the complement of $A_v$. Then
\begin{eqnarray}\label{avavc}
\frac{\Pr(A_v^c)}{\Pr(A_v)} = \prod_{\omega \in Y_v} \left[1 + \displaystyle{ \frac{\sum_{S \subseteq \omega} (1 - f_S^+/f_S) \cdot  \beta_S f_S } { W^+_\omega}}\right] \notag.
\end{eqnarray}

{\it Proof.} Since $H$ is linear, the sets in $Y_v$ are disjoint, and since $H$ is locally sparse, each edge of $H$ contains at most one vertex from at most one set in $Y_v$. 
Therefore, conditioning on $I_v$, we may add any subset of vertices in $\bigcup_{\omega \in Y_v}  \omega$ to $I$ and we may add $v$ to $I$ iff none of the $\omega \in Y_v$ with $|\omega| = r$ is completly added to $I$. Therefore by summing over all subsets of $\bigcup_{\omega \in Y_v}$ to add to $I$ and  by definition (\ref{measure}),
\begin{align}\label{avequation}
	\Pr(A_v) &=   \frac{ \prod_{\omega \in Y_v}     \sum_{ S \subseteq \omega, |S| < |r| }   \beta_S f_S^+ } { \prod_{\omega \in Y_v} \sum_{S \subseteq \omega, |S| < r } \beta_Sf_S^+   + \sum_{S \subseteq \omega } \beta_S f_S }\\
& = \frac{ \prod_{\omega \in Y_v}     \sum_{ S \subseteq \omega }   \beta_S f_S^+ } { \prod_{\omega \in Y_v} \sum_{S \subseteq \omega } \beta_Sf_S^+   + \sum_{S \subseteq \omega } \beta_S f_S }\\
 &=   \frac{\prod_{\omega \in Y_v}  W^+_\omega }{\prod_{\omega \in Y_v} W^+_\omega  \; \; + \; \;  \prod_{\omega \in Y_v} W_\omega }.
\end{align}
The equality of the first and second line follows from the fact that $\beta_{r + 1} = 0$ and so $f_S^+ = 0$ when $|S| = r$.  From the left side in (\ref{avavc}), we obtain
\begin{align*}
\prod_{\omega \in Y_v} \left[1 + \displaystyle{ \frac{\sum_{S \subseteq \omega} (1 - f_S^+/f_S) \cdot \beta_S f_S } { W^+_\omega}}\right]  = \prod_{\omega \in Y_v} \left[\displaystyle{  \frac{ W_{\omega}^+ +\sum_{S \subseteq \omega} (\beta_S f_S - \beta_S f_S^+) } { W^+_\omega}}\right]  = \prod_{\omega \in Y_v} \frac{\displaystyle{ W_\omega }}{ { W_\omega^+}}.
\end{align*}
By (\ref{avequation}), this is precisely $\Pr(A_v^c)/\Pr(A_v)$, as required. $\bs$

\bigskip
\medskip

{\bf Claim 4.}  For $v \in V(H)$ and $\omega \in Y_v$,
\begin{equation*}
   \frac{W_\omega}{W^+_\omega} \leq 1 + \frac{1}{2^r -1} + o_d(1).
\end{equation*}
{\it Proof.} Recall from (\ref{parameters}) that $\beta_{r + 1} = 0$ and $\beta_i \rightarrow 1$ as $d \rightarrow \infty$ for $0 \leq i \leq r$. Therefore if $|S| < r$, $f_S \rightarrow 1 = f_S^+$ as $d \rightarrow \infty$, and if $|S| = r$, $f^+_S = 0$ and $f_S \rightarrow 1$ as $d \rightarrow \infty$. So, if $|\omega| < r$, then 
$W_{\omega}/W_{\omega}^+ \rightarrow 1$ as $d \rightarrow \infty$. On the other hand, if $|\omega| = r$, then
\begin{align*}
\limsup_{d \rightarrow \infty}  \frac{W_\omega}{W^+_\omega} & \leq \limsup_{d \rightarrow \infty} \frac{  f_{\omega} \beta_{\omega}  + \sum_{S \subsetneq \omega} f_S \beta_S}{ \sum_{S \subsetneq \omega} f^+_S \beta_S}\\
&= \limsup_{d \rightarrow \infty}  \frac{\beta_{\omega} + \sum_{S \subsetneq \omega} \beta_S}{\sum_{S \subsetneq \omega} \beta_S} \\
&= 1 + \limsup_{d \rightarrow \infty}  \frac{1}{\sum_{S \subsetneq \omega} \beta_S/\beta_{\omega}} \; \; \leq \; \;  1 + \frac{1}{2^r - 1}.
\end{align*}
Here we used the fact that $\beta_{\omega} \leq \beta_S$ for $S \subseteq \omega$. $\bs$

\bigskip
\medskip

{\bf Claim 5.} For $v \in V(H)$ and $I \in \mathcal{I}(H)$,
\begin{equation}\label{xviv}
	\mathrm{E}(X_v|I_v) \geq \Pr(A_v) \cdot d^{\frac{1}{r}} +  (1 - o_d(1))\frac{(2^r - 1)}{r 2^r} \cdot  \Pr(A_v^c) \log \frac{\Pr(A_v^c)}{\Pr(A_v)}.
\end{equation}
In particular, for all $v \in V(H)$, 
\begin{equation}\label{xvlower}
\mathrm{E}(X_v) \geq (1  - o_d(1))  \frac{2^r - 1}{r^2 2^r} \cdot \log d.
\end{equation}

{\it Proof.} Since $-\log(\beta_{j + 1}/\beta_j) \leq -\log \beta_{j + 1}$ and $1 - x \leq -\log x$ for $x > 0$, for $0 < |S| < r$,
\begin{eqnarray*}
1 - \frac{f_S^+}{f_S} = 1 - \prod_{j = 0}^{|S|} \Bigl(\frac{\beta_{j + 1}}{\beta_j}\Bigr)^{{|S| \choose j}} \leq \sum_{j = 0}^{|S|} {|S| \choose j} (r - j) \delta^{r - j} \leq r \delta^{r - |S|}  \sum_{j = 0}^{|S|}  { |S| \choose j } \delta^{ |S| - j }
\end{eqnarray*}
By the binomial theorem, this is exactly
\begin{eqnarray*}
 r(1 + \delta)^{|S|} \delta^{r - |S|} \leq r(1 + \delta)^r \cdot \delta^{r - |S|}.
\end{eqnarray*}
For $|S| = r$, we have $1 - f_S^+/f_S \leq 1 \leq r(1 + \delta)^r \delta^{r - |S|}$, and for $|S| = 0$, $1 - f_S^+/f_S = 0$. Using $\log(1 + t) \leq t$ for $t > -1$, Claim 3 gives
\begin{eqnarray}\label{fromclaim3}
\log \frac{\Pr(A_v^c)}{\Pr(A_v)} \leq r(1 + \delta)^r  \cdot \sum_{\omega \in Y_v}\frac{1}{W_{\omega}^+} \sum_{\emptyset \subsetneq S \subseteq \omega}  \delta^{r-|S|}\beta_S f_S.
\end{eqnarray}
By definition (\ref{xv}) of $X_v$ 
\begin{align*}
	\mathrm{E}(X_v|I_v) &= \Pr(A_v) \cdot d^{\frac{1}{r}} +  \Pr(A_v^c) \cdot \sum_{\omega \in Y_v} \sum_{0 \subseteq S \subseteq \omega} P( I \cap \omega = S) \sum_{0 \subsetneq T \subseteq S} \delta^{r - |T|}  \\ 
& = \Pr(A_v) \cdot d^{\frac{1}{r}} +  \Pr(A_v^c) \cdot \sum_{\omega \in Y_v} \frac{1}{W_\omega} \sum_{0 \subsetneq S \subseteq \omega}
 \beta_ S f_S \sum_{0 \subsetneq T \subseteq S} \delta^{r - |T|}  
\end{align*}
We now apply (\ref{fromclaim3}) and Claim 4 to this expression to obtain 
\begin{align*}
\mathrm{E}(X_v|I_v) & \geq \Pr(A_v) \cdot d^{\frac{1}{r}} +  (1 - o_d(1)) \cdot \frac{\Pr(A_v^c)}{1 + 1/(2^r - 1)}  \cdot \sum_{\omega \in Y_v} \frac{1}{W_\omega^+} \sum_{0 \subsetneq S \subseteq \omega}
 \beta_ S f_S  \delta^{r - |S|}  \\ 
&\geq  \Pr(A_v) \cdot d^{\frac{1}{r}} +  (1 - o_d(1)) \cdot \frac{\Pr(A_v^c)}{(1 + 1/(2^r - 1))r(1 + \delta)^r} \cdot \log \frac{\Pr(A_v^c)}{\Pr(A_v)} \\[0.7em]
&= \Pr(A_v) \cdot d^{\frac{1}{r}} +  (1 - o_d(1)) \cdot \frac{2^r - 1}{r2^r(1 + \delta)^r} \cdot \Pr(A_v^c) \log \frac{\Pr(A_v^c)}{\Pr(A_v)}.
\end{align*}
The fact that $(1 + \delta)^r = 1 + o_d(1)$ finishes the proof of (\ref{xviv}). To obtain (\ref{xvlower}), note that the lower bound in (\ref{xviv}) 
is minimized as a function of $\Pr(A_v)$ when
\[
\Pr(A_v) = d^{-\frac{1}{r} + o_d(1)}
\]
in which case the minimum for each $v \in V$ is at least
\[ (1  - o_d(1))  \frac{2^r - 1}{r^2 2^r} \cdot \log d.\]
We conclude $\mathrm{E}(X_v|I_v)$ is at least this value for all $v \in V(H)$, which proves (\ref{xvlower}). $\bs$

\subsection{Proof of Theorem \ref{main}}

To prove Theorem \ref{main}, it is enough to show
\[ \alpha(H) \geq (1 - o_d(1)) \cdot \frac{(r - 1)(2^r - 1)}{r^2 2^r} \cdot \delta n\]
due to the definition of $\delta$. The upper bound in Claim 2 and the lower bound (\ref{xviv}) in Claim 5 give
\[ (1  - o_d(1)) \cdot \frac{2^r - 1}{r^2 2^r} \cdot n\log d  \; \; \leq \; \; \mathrm{E}(X) \; \; \leq \; \; (1 + o_d(1)) \frac{\log d}{\delta (r-1)} \cdot \alpha(H).\]
Rearranging terms, we obtain
\begin{eqnarray*}
\alpha(H) &\geq& (1 - o_d(1))  \cdot \frac{(r - 1)(2^r - 1)}{r^2 2^r} \cdot \delta n,
\end{eqnarray*}
as required.  \qed

\section{Concluding remarks}

$\bullet$ By a suitable change of constants, Theorem \ref{main} can be applied to the case of $(r + 1)$-uniform hypergraphs of {\em average} degree $d$, by simply deleting all high degree vertices.

\bigskip

$\bullet$ If $H$ is a random $d$-regular $(r + 1)$-uniform $n$-vertex hypergraph, where $r \geq 1$ and $d \geq 1$, then it may be shown that with positive probability, $H$ is locally sparse if $d$ is not too large relative to $n$ and, with high probability, 
\begin{equation}\label{randomreg}
\alpha(H) = (1 + o_d(1))(1 + 1/r)^{1/r} \cdot \Bigl(\frac{\log d}{d}\Bigr)^{\frac{1}{r}} \cdot n.
\end{equation}
It is a major open question as to whether this is tight in the case of graphs, namely, does 
every locally sparse (i.e. triangle-free) $n$-vertex graph of maximum degree $d$ have an independent set of size at least $(2 + o_d(1))n(\log d)/d$ as $d \rightarrow \infty$? Shearer~\cite{Sh} proved 
that every triangle-free $n$-vertex graph $G$ of maximum degree $d$ has $\alpha(G) \geq n(d\log d - d + 1)/(d - 1)^2$ for $d \geq 2$, which is a factor two 
away from (\ref{randomreg}). We make the following more modest conjecture: 
\begin{conj} \label{mainconj}
For $r \geq 2$, there exists a constant $a_{r} > 0$ such that $a_r = 1 - o_r(1)$ as $r \rightarrow \infty$ and every 
$(r + 1)$-uniform $n$-vertex linear hypergraph $H$ with maximum degree $d \geq 1$ has
\begin{equation}\label{alphaguess}
	\alpha(H) \geq (a_r - o_d(1)) \cdot \Bigl(\frac{\log d}{d}\Bigr)^{\frac{1}{r}} \cdot n \quad \quad \mbox{ as } d \rightarrow \infty.
 \end{equation}
\end{conj}
It appears to be difficult to determine whether this is true with $a_r = (1 + 1/r)^{1/r}$ for any $r \geq 1$, or whether one may take $a_r > 1$ for any $r \geq 1$.
For $d$-regular hypergraphs of ``large girth'', one obtains $a_r = r^{-1/r}$, as shown in~\cite{NV}, which shows Conjecture \ref{randomreg} is true for this class of hypergraphs.

\bigskip 

$\bullet$  The following strengthening of Conjecture \ref{mainconj} was proposed for linear triple systems (see Conjecture 1.2 in~\cite{EV} and see~\cite{KMV}): 
\begin{conj} \label{steiner}
If $H$ is an $n$-vertex linear triple system, then 
\begin{equation}
	\alpha(H) \geq (1 - o_n(1)) \cdot \sqrt{3n\log n} \quad \mbox{ as }n \rightarrow \infty.
 \end{equation}
\end{conj}
In particular, this states that every Steiner triple system has independence number asymptotically at least as large as the expected independence number of random triple system of the same density on the same set of vertices -- see the third concluding remark above. Note that any linear $n$-vertex triple system has maximum degree at most $d = (n - 1)/2$ and, as in (\ref{randomreg}), a random $d$-regular triple system has expected independence number $(1 + o_n(1))\sqrt{3n\log n}$ as $n \rightarrow \infty$.  

\bigskip

$\bullet$ As remarked in Section \ref{intuition}, if $I$ is a uniformly selected independent set in a locally sparse or linear $(r + 1)$-uniform $n$-vertex 
hypergraph $H$ of maximum degree $d$, then we lack a good upper bound for $\mathrm{E}(|\partial_rI|)$. We make the following conjecture, which may be 
independently interesting:

\begin{conj}\label{uniform}
For $r \geq 2$ there exists $C_r > 0$ such that if $H$ is any $n$-vertex linear $(r + 1)$-uniform hypergraph of maximum degree $d \geq 1$ and $I \in \mathcal{I}(H)$ is sampled uniformly, then
\begin{equation}\label{guess}
\mathrm{E}(|\partial_r I|) \leq C_r \cdot d \alpha(H)^r n^{1-r}.
\end{equation}
\end{conj}

Note that if $I$ is a uniformly sampled subset of $V(H)$ of size $\alpha(H)$, then $\mathrm{E}(|\partial I|) \leq d\alpha(H)^r n^{1 - r}$, so 
the conjecture is proposing that the expected value of $|\partial_r I|$ is not much larger than if $I$ were any subset of vertices of $H$  
of size $\alpha(H)$ chosen uniformly at random. In that case, following Section \ref{intuition}, we would have for some $b_r > 0$:
\[ b_r \cdot n\log d  \leq \mathrm{E}(X) \leq d^{1/r} \cdot \mathrm{E}(|I|) + C_r \cdot d \alpha(H)^r n^{1 - r}   \]
and therefore for some $B_r > 0$ where $B_r \approx (b_r/C_r)^{1/r}$, we have 
\[ \alpha(H) \geq B_r \cdot \frac{n(\log d)^{1/r}}{d^{1/r}} \]
which gives a bound of the same order of magnitude as in Theorem \ref{main}. If $C_r$ is sub-exponential in $r$ -- and perhaps $C_r = 1 + o_d(1)$ -- then 
this shows Conjecture \ref{uniform} implies Conjecture \ref{mainconj}. 
If $H$ is an $n$-vertex Steiner triple system, then $d = (n - 1)/2$ and for every set $I \subseteq V(H)$, $\partial_2 I = {I \choose 2}$ and so $|\partial_2 I| \leq d\alpha(H)^2/n$, consistent with (\ref{guess}). This was used in~\cite{KMV} to give good lower bounds on the independence number 
of Steiner systems.

\medskip

$\bullet$ Duke, Lefmann and R\"{o}dl~\cite{DLR} showed that by randomly sampling vertices, one may easily obtain a lower bound on $\alpha(H)$ of the same order of magnitude
as in Theorem \ref{main} for all linear $n$-vertex $(r + 1)$-uniform hypergraphs $H$ of maximum degree $d \geq 1$ when $r \geq 2$. Specifically, 
if we sample vertices with probability $p$, then the expected number of $k$-cycles whose vertices are all sampled 
is at most $p^{rk} ((r + 1)d)^{k - 1} n$, whereas the expected number of edges of sampled vertices on a sampled vertex is at most $p^{r}d$.  
Choosing $p = (2(r + 1)d)^{-2/(3r - 1)}$ and using probabilistic methods, we obtain with positive probability an induced 
subgraph $H'$ with $(1 - o_d(1))pn$ vertices and in which all vertices have degree at most roughly $D = (1 + o_d(1))p^r d$. Applying Theorem \ref{main} 
to $H'$, we obtain 
\[ \alpha(H) \geq \alpha(H') \geq c_r \cdot \Bigl(\frac{\log D}{D}\Bigr)^{\frac{1}{r}} \cdot (1 - o_d(1))pn  \geq (1 - o_d(1)) \cdot \Bigl(1 - \frac{2r}{3r-1}\Bigr)^{\frac{1}{r}}  \cdot c_r \cdot \Bigl(\frac{\log d}{d}\Bigr)^{\frac{1}{r}} \cdot n.\]
In particular, asymptotically as $d,r \rightarrow \infty$, this gives the same lower bound as in Theorem \ref{main} for linear hypergraphs.

\end{document}